\documentclass[a4paper]{article}

\usepackage{graphicx, here}
\usepackage{fancyheadings}
\usepackage{fonte}
\usepackage{timesmt}
\usepackage[usenames,dvipsnames,table]{xcolor}
\usepackage{amsmath, amssymb, amsfonts, amsmath, bm}

\newcommand{\X}{\mathbb{X}}
\newcommand{\R}{\mathbb{R}}

\newcommand \m {\dot{m}}
\newcommand \mr {\dot{m}_{r}}
\newcommand{\rpc}{r_{2p}}
\newcommand{\rtc}{r_{2t}}
\newcommand{\rc}{r_{3}}
\newcommand{\rt}{r_{4}}
\newcommand{\rpt}{r_{5p}}
\newcommand{\rtt}{r_{5t}}
\newcommand{\bc}{b_3}
\newcommand{\bt}{b_4}
\newcommand{\betac}{\beta_3}
\newcommand{\alphat}{\alpha_4}
\newcommand{\Lxa}{L_{x1}}
\newcommand{\Lya}{L_{y1}}
\newcommand{\Lza}{L_{z1}}
\newcommand{\Lxb}{L_{x2}}
\newcommand{\Lyb}{L_{y2}}
\newcommand{\Lzb}{{L_{z2}}}
\newcommand{\dphx}{\Delta P_{HX}}

\begin{document}

\thispagestyle{fancy}
\lhead{{\ninerm \bf  EngOpt 2016 - $5^{th}$ International Conference on Engineering Optimization} \\
{\eightrm Iguassu Falls, Brazil, 19-23 June 2016.}}
\begin{center}
{\large {\elevenrm \bf Design of a commercial aircraft environment control system using Bayesian optimization techniques}}\\~\\
{\tenrm \bf Paul Feliot\textsuperscript{1,2}, Yves Le Guennec\textsuperscript{1}, Julien Bect\textsuperscript{2,1} and Emmanuel Vazquez\textsuperscript{2,1}}
\\~\\
{\ninerm
\textsuperscript{1} Institut de Recherche Technologique SystemX, 8 avenue de la Vauve, Palaiseau, France (firstname.surname@irt-systemx.fr)\\
\textsuperscript{2} Laboratoire des Signaux et Syst{\`e}mes (L2S), CentraleSup{\'e}lec, CNRS, Univ. Paris-Sud, Universit{\'e} Paris-Saclay, Gif-sur-Yvette, France (firstname.surname@centralesupelec.fr)
}
\end{center}

\begin{abstract}
In this paper, we present the application of a recently developed algorithm for Bayesian multi-objective optimization to the design of a commercial aircraft environment control system (ECS). In our model, the ECS is composed of two cross-flow heat exchangers, a centrifugal compressor and a radial turbine, the geometries of which are simultaneously optimized to achieve minimal weight and entropy generation of the system. While both objectives impact the overall performance of the aircraft, they are shown to be antagonistic and a set of trade-off design solutions is identified. The algorithm used for optimizing the system implements a Bayesian approach to the multi-objective optimization problem in the presence of non-linear constraints and the emphasis is on conducting the optimization using a limited number of system simulations. Noteworthy features of this particular application include a non-hypercubic design domain and the presence of hidden constraints 
due to simulation failures. \\

\end{abstract}

\section{Introduction}
\label{sec:intro}

The purpose of the environment control system (ECS) of a commercial aircraft is to provide a certain level of comfort to the passengers by regulating the temperature and pressure of the air injected into the cabin. The system is based on an inverse Brayton thermodynamic cycle. Hot and pressurized air is taken from the engines at the compressor stage through the bleed and ram air from the outside of the aircraft is used as coolant. For safety reasons, the hot air from the engines passes through a first heat exchanger where it is cooled down below the critical fuel ignition temperature. Then it is pressurized through a compressor and cooled again using a second heat exchanger. It then passes through a turbine where work is extracted to propel the compressor. The cooled and expanded air exiting the turbine is eventually mixed with hot air from the first heat exchanger outflow to reach the desired temperature and pressure before injection into the cabin.

The design of an optimal ECS is a complex problem in practice. It has been addressed in previous studies under different optimality conditions and modelling assumptions (see, e.g., \cite{vargas2001thermodynamic,bejan2001need,perez2002optimization}). In their article, P{\'e}rez-Grande and Leo \cite{perez2002optimization} study an aircraft-on-cruse scenario and propose a one dimensional model of the two heat exchangers. The system is designed in order to achieve minimal mass and entropy generation, two objectives that are shown to be antagonistic and which both affect the overall performance of the aircraft.

In this article, we extend their work by considering also the sizing of the rotating machines and by considering an aircraft-on-ground scenario, which corresponds to the most critical situation for the ECS in terms of cold production, and is therefore dimensioning. The design optimization of the system is performed using the BMOO algorithm \cite{feliot2016bayesian}, which implements a Bayesian approach to the multi-objective optimization problem in the presence of non-linear constraints. The problem consists in finding an approximation of the set
\begin{equation*}
\label{eq:Gamma}
  \Gamma = \{x \in \X : c(x) \leq 0 \text{ and } \nexists\, x'\in\X
  \text{\ such that\ } c(x') \leq 0 \text{ and } f(x') \prec f(x) \}  
\end{equation*}
where $\X \subset \R^d$ is the search domain, $c = (c_i)_{1 \leq i \leq q}$ is a
vector of constraint functions ($c_i:\X\to\R$), $c(x) \leq 0$ means that
$c_{i}(x) \leq 0$ for all $1\leq i \leq q$, $f = (f_j)_{1 \leq j \leq p}$ is a
vector of objective functions to be minimized ($f_j:\X\to\R$), and $\prec$
denotes the Pareto domination rule (see, e.g., \cite{fonseca1998multiobjective}).

The article is organized as follows. First we detail in Section~\ref{sec:analysis} the model that is used to estimate the performances and main characteristics of the system. A one-dimensional analysis is performed to establish the state equations of the system and link the physical values of interest to the geometrical parameters of the system components. Then we introduce in Section~\ref{sec:optim} the optimization algorithm that is used to conduct the system optimization. The results of the optimization are analyzed and possible directions for future work are discussed. Finally, conclusions are drawn in Section~\ref{sec:conclusion}.

\section{Thermodynamic analysis of the ECS}
\label{sec:analysis}

\subsection{Sizing scenario}
\label{sec:scenario}

\begin{figure*}
\centerline{
\includegraphics[width=11cm]{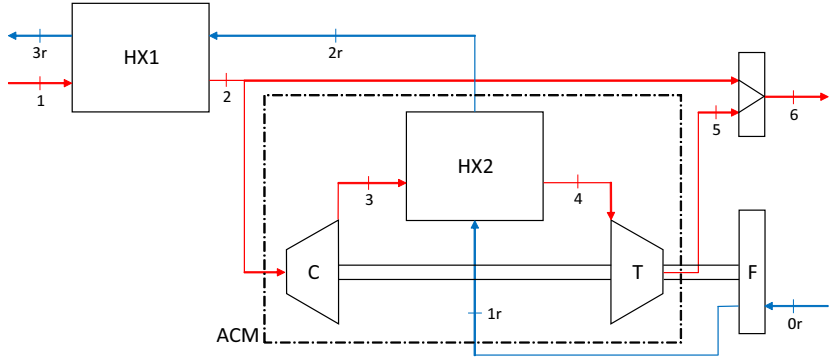}
}
\caption{Architecture of the environment control system of a commercial aircraft}
\label{fig:ecs}
\end{figure*}

The architecture of the ECS is represented on Figure~\ref{fig:ecs}. Bleed air from the engines arrives into the system at location 1. Ram air from the outside of the aircraft is levied at location 0r and is used as coolant. The hot air enters a first heat exchanger where it is cooled down below the fuel ignition temperature. A by-pass at location 2 then permits to regulate the system by controlling the air flowrate entering the air cycle machine (ACM). Cooled and expanded air exits the ACM at location 5 and is mixed with warm air from the by-pass to reach the desired pressure and temperature before injection into the cabin at location 6.

In practice, the system must be able to satisfy strict specifications under different environmental conditions and operating situations. In this work, we consider a scenario where the aircraft is on ground, full of passengers, equipments running, and with an outside temperature of $50^\circ$C. In that situation, the ECS must be able to maintain the cabin temperature at $T_c = 24^\circ$C. This scenario corresponds to the most demanding specification in terms of cold production, and is therefore dimensioning for the system. Formally, this means that the ECS must be able to dissipate enough heat to compensate for the thermal power $\mathcal{P}_{HT}$ produced by the passengers, the crew, the equipments and the environment:
\begin{equation*}
\mathcal{P}_{HT} = \mathcal{P}_{out} + \mathcal{P}_{eq} + N_{pax}\mathcal{P}_{pax} + N_{crew}\mathcal{P}_{crew},
\end{equation*} 
where $\mathcal{P}_{out}$ is the outside flow dissipation, $\mathcal{P}_{eq}$ is the thermal power produced by the equipments, $\mathcal{P}_{pax}$ and $\mathcal{P}_{crew}$ are the thermal powers produced by a passenger and by a crewmember and $N_{pax}$ and $N_{crew}$ are the number of passengers and the number of crewmembers in the aircraft.

In this scenario, the by-pass is wide open so that all the air from location 2 goes to the ACM. Also, there is no relative velocity between the aircraft and the ambient air when it is grounded and therefore there is no natural coolant flowrate. In this work, we consider a system where the ram flowrate is created by an auxiliary fan placed at the ram air entrance and powered by the turbine of the ACM. The sizing of this auxiliary fan is not taken into account in this study but we will assume that the ram flowrate can be controlled.

\subsection{Heat exchangers}

We now detail the model that is used to emulate the system. For the heat exchangers HX1 and HX2, we use a model from \cite{perez2002optimization}. The two heat exchangers are compact cross-flow heat exchangers with unmixed fluids. For this kind of heat exchangers, the energy exchanged per unit time between the ram and bleed air can be formulated as:
\begin{eqnarray}
\m c_p (T_{t1} - T_{t2}) &=& \mr c_p (T_{t3r} - T_{t2r}), \\
\m c_p (T_{t3} - T_{t4}) &=& \mr c_p (T_{t2r} - T_{t1r}),
\end{eqnarray}
where $\m$ and $\mr$ denote respectively the bleed and ram air flowrates, $c_p$ is the thermal capacity of the air and is assumed constant, and $T_{ti}$ represents the stagnation temperature at location $i \in \{ 1,2,3,4,5,1r,2r,3r \}$. Besides, we can define the efficiencies $\epsilon_1$ and $\epsilon_2$ of the two heat exchangers HX1 and HX2 as the ratio between the energy effectively exchanged and the total energy exchangeable: 
\begin{eqnarray}
c_p (T_{t1} - T_{t2}) &=& \epsilon_1 c_p (T_{t1} - T_{t2r}), \\
c_p (T_{t3} - T_{t4}) &=& \epsilon_2 c_p (T_{t3} - T_{t1r}).
\end{eqnarray}

Note that Eq.(3-4) only hold when $\m \leq \mr$. Otherwise, $T_{t1} - T_{t2r}$ should be replace by $T_{t3r} - T_{t2r}$ in Eq.(3) and $T_{t3} - T_{t1r}$ should be replaced by $T_{t2r} - T_{t1r}$ in Eq.(4). The efficiencies $\epsilon_1$ and $\epsilon_2$ of the two heat exchangers depend on their geometry and we use the $\epsilon$-Ntu model detailed in \cite{perez2002optimization} to estimate them.

In this study, the pressure drops as the air passes through the heat exchangers are considered constant: 
\begin{eqnarray}
P_{t2} - P_{t1} &=& \dphx , \\
P_{t4} - P_{t3} &=& \dphx ,
\end{eqnarray}  
where $P_{ti}$ represents the stagnation pressure at location $i$ and $\dphx$ is a constant pressure loss. While this permits to simplify the model, it is inaccurate because the pressure losses do also depend on the geometries of the heat exchangers and are responsible for a non-negligible proportion of the entropy generated by the system. In particular, the frictions are expected to rise as the volume of the heat exchangers decreases, which further increases the necessary balance between entropy generation and mass. Again, the reader is referred to \cite{perez2002optimization} for a discussion on a possible model of the pressure drops.

\subsection{Compressor and turbine}

\begin{figure}
\centerline{
\includegraphics[height=4cm]{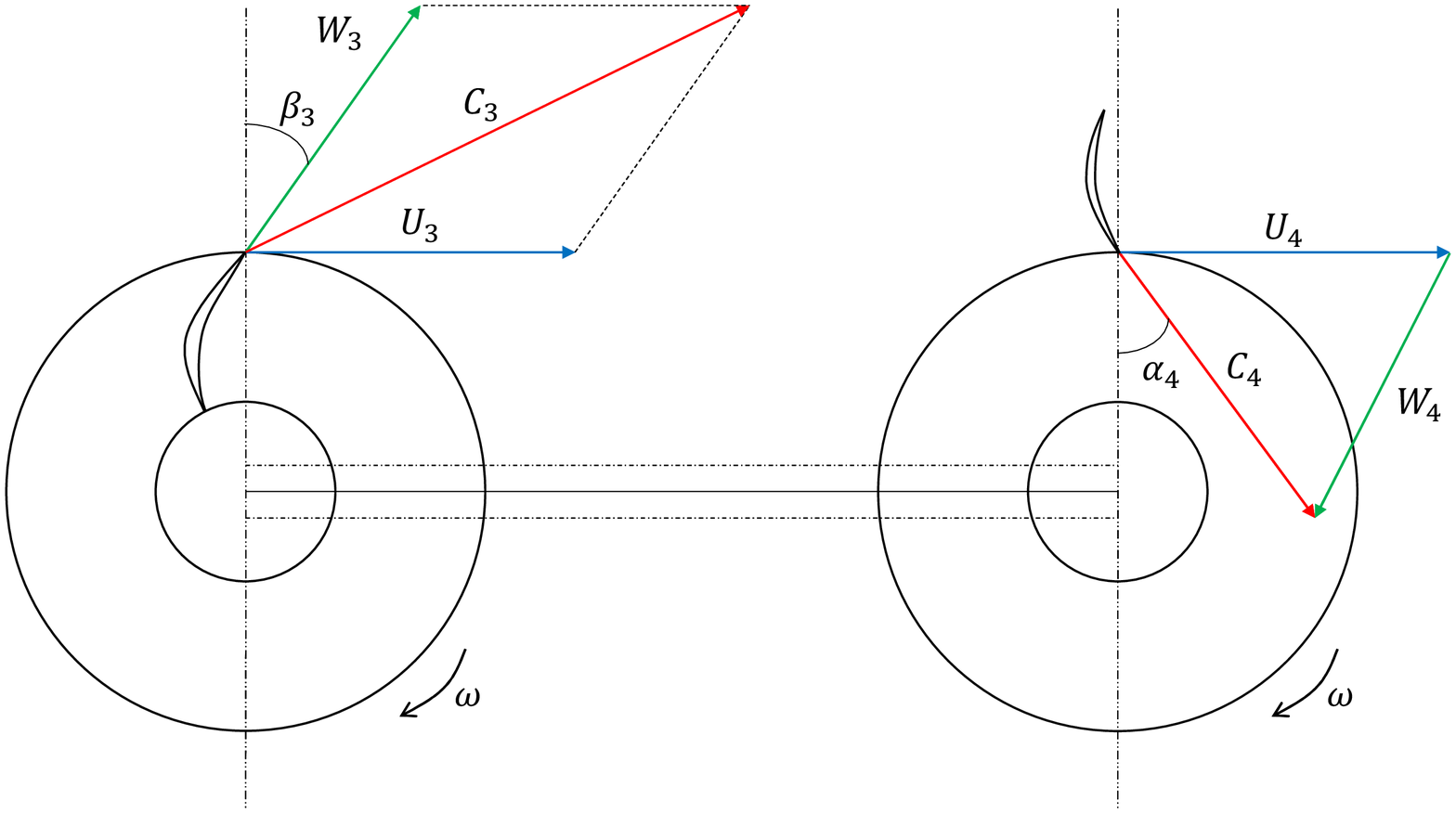}
}
\caption{Compressor (left) and turbine (right) velocity triangles.}
\label{fig:triangles}
\end{figure}

We consider a centrifugal compressor and an axial flow turbine. It is assumed that the air enters the compressor axially and exits parallel to the blades (the slip factor is neglected), with an angle $\betac$ as illustrated on Figure~\ref{fig:triangles}. Similarly, it is assumed that the air enters the turbine with an angle $\alphat$ corresponding to the stator blades  angle and goes out axially. The rotational speed of the shaft linking the compressor with the turbine is denoted $\omega$. The letters $C$, $U$ and $W$ on Figure~\ref{fig:triangles} denote respectively the air absolute velocity vector, the blade tip speed vector and the relative velocity vector, such that $C = U + W$. In the following, the subscripts $u$, $m$, and $x$ will stand respectively for the tangential, meridional and axial components of the velocity vectors.

The power exchanged per unit time between the machines and the fluid and the change of momentum of the fluid are related by Euler's theorem as:
\begin{eqnarray*}
\dot{W}_C &=& \m(U_3 C_{3u} - U_2 C_{2u}) , \\
\dot{W}_T &=& \m(U_5C_{5u} - U_4C_{4u}) ,
\end{eqnarray*}
where $\dot{W}_C$ and $\dot{W}_T$ denote respectively the power received by the fluid from the compressor and from the turbine. Note that with this convention $\dot{W}_T < 0$ and $\dot{W}_C > 0$. The turbine is converting part of the fluid energy into rotation speed, and the compressor is augmenting the energy of the fluid through its mechanical work, which increases the temperature of the fluid. 

Under the assumption that the flow enters the compressor and exits the turbine axially, the tangential components of the air velocity vectors are neglected and $C_{2u} = C_{5u} = 0$. The Euler theorem simplifies to:
\begin{eqnarray*}
\dot{W}_C &=& \m U_3 C_{3u} , \\
\dot{W}_T &=& -\m U_4 C_{4u} ,
\end{eqnarray*}

The conservation of the mass at locations 3 and 4 gives the following additional equations, relating the air flowrate to its velocity and a control surface.
\begin{eqnarray*}
\m &=& 2 \pi \rho r_3 b_3 W_{3m} , \\
\m &=& 2 \pi \rho r_4 b_4 C_{4m} ,
\end{eqnarray*}
where $\rho$ is the air density (which is assumed constant), $r_3$ and $r_4$ are respectively the compressor outlet blade radius and turbine inlet blade radius, and $b_3$ and $b_4$ are respectively the compressor and turbine tip blade heights (see Figure~\ref{fig:parameters}). Using the velocity triangles of Figure~\ref{fig:triangles}, the Euler theorem can then be rewritten using the rotating machines geometries and rotational speed:
\begin{eqnarray}
\dot{W}_C &=& \m \left( r_3^2\omega^2 - \frac{\m\tan(\beta_{3})}{2\pi \rho b_3}\omega \right) , \\
\dot{W}_T&=& -\frac{\m^2\tan(\alpha_4)}{2\pi \rho b_4}\omega .
\end{eqnarray}

Besides, the work extracted from the turbine is used to propel the compressor and the auxiliary fan. Writing down the conservation of energy per unit time we get:
\begin{equation}
\dot{W}_C + \dot{W}_T + \frac{1}{\eta_F}\frac{\mr^3}{2 \rho^2 A_r^2} = 0 ,
\label{eq:omega-poly}
\end{equation}
where $A_r$ is a control surface at the ram air entrance and $\eta_F$ is the ratio between the kinetic energy per unit time produced by the auxiliary fan and the power furnished by the turbine to the fan.

The powers $\dot{W}_C$ and $\dot{W}_T$ can also be expressed as functions of the stagnation temperatures by considering the change in total enthalpy of the fluid passing through the rotating machines (the other contributions are neglected):
\begin{eqnarray}
\dot{W}_C &=& \eta_C \m c_p(T_{t3}-T_{t2}) , \\
\dot{W}_T &=& \frac{1}{\eta_T} \dot{m} c_p(T_{t5}-T_{t4}) ,
\end{eqnarray} 
where $\eta_C$ and $\eta_T$ are respectively the compressor and turbine isentropic efficiencies. Finally, the stagnation pressure ratios are given by the isentropic relations:
\begin{eqnarray}
\frac{P_{t3}}{P_{t2}} &=& \left(1 + \eta_C\frac{T_{t3}-T_{t2}}{T_{t2}}\right)^{\frac{\gamma}{\gamma-1}} , \\
\frac{P_{t5}}{P_{t4}} &=& \left(1 + \frac{1}{\eta_T}\frac{T_{t5}-T_{t4}}{T_{t4}}\right)^{\frac{\gamma}{\gamma-1}} ,
\end{eqnarray}
where $\gamma$ is the air isentropic coefficient.

\subsection{Mass and entropy generation rate of the system}

\begin{figure}
\centerline{
\includegraphics[height=4cm]{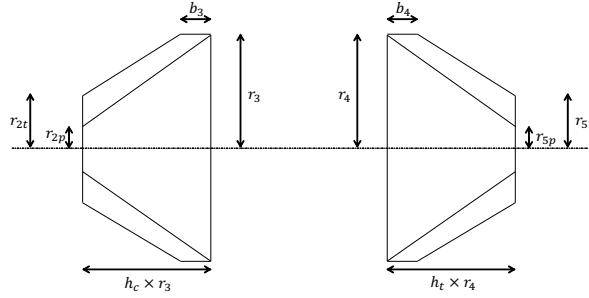}
}
\caption{Geometrical parametrization of the compressor (left) and of the turbine (right).}
\label{fig:parameters}
\end{figure}

The mass of the system can be approximated by considering estimates of the volumes of its components and representative densities. For the two heat exchangers, we consider rectangular volumes and a representative density $\rho_{HX}$ (see Table~\ref{tab:static}).
\begin{eqnarray*}
\mathcal{M}_{HX1} &=& \rho_{HX} \cdot \Lxa \Lya \Lza, \\
\mathcal{M}_{HX2} &=& \rho_{HX} \cdot \Lxb \Lyb \Lzb, 
\end{eqnarray*}
where $\Lxa$, $\Lya$, $\Lza$, $\Lxb$, $\Lyb$, $\Lzb$ are the heat exchangers dimensions (see Table~\ref{tab:variables}), and  $\mathcal{M}_{HX1}$ and $\mathcal{M}_{HX2}$ are the mass of the heat exchangers. For the compressor and turbine of the ACM, we consider separately the volumes of the blades and the volume of the machine body (see Figure~\ref{fig:parameters}):
\begin{eqnarray*}
V_{C, blade} &=& e_c \left( \frac{h_c \rc(\rc-\rpc)}{2} - \frac{(\rc-\rtc)(h_c\rc-\bc)}{2} \right), \\
V_{C, body} &=& \frac{\pi \rc^2 h_c(\rc+\rpc)}{3} - \frac{\pi h_c \rpc^3}{3},
\end{eqnarray*}
for the compressor, and
\begin{eqnarray*}
V_{T, blade} &=& e_t \left( \frac{h_t \rt(\rt-\rpt)}{2} - \frac{(\rt-\rtt)(h_t\rt-\bt)}{2} \right), \\
V_{T, body} &=& \frac{\pi \rt^2 h_t(\rt+\rpt)}{3} - \frac{\pi h_t \rpt^3}{3},
\end{eqnarray*}
for the turbine, where $e_c$ and $e_t$ are the compressor and turbine blades thickness, and $h_c$ and $h_t$ are aspect ratios. The mass of the system is then given by the following.
\begin{eqnarray}
\mathcal{M} &=& \mathcal{M}_{HX1} + \mathcal{M}_{HX2} + \rho_{steel}(Z_C V_{C, blade} + V_{C, body}) + \rho_{steel}(Z_T V_{T, blade} + V_{T, body}),
\end{eqnarray}
where $Z_C$ and $Z_T$ are respectively the number of blades of the compressor and of the turbine.

The entropy generation rate of the ECS is the sum of the contributions along the bleed stream and along the ram stream, from entrance to exit:
\begin{eqnarray}
\dot{\mathcal{S}} &=& \m\left( c_p\log \frac{T_5}{T_a} - R \log \frac{P_5}{P_a} \right) + \mr\left( c_p\log \frac{T_{3r}}{T_{a}} - R \log \frac{P_{3r}}{P_{a}} \right),
\end{eqnarray}
where $R$ is the perfect gas constant, $T_a$ and $P_a$ are the ambient temperature and pressure, and $T_i$ and $P_i$ are respectively the static temperature and the static pressure at location $i \in \{ 5,3r \}$. The equations giving the static properties for the bleed stream are gathered in Table~\ref{tab:static} in the additional material. For the ram stream, it is assumed that the Mach number remains low through the heat exchangers. Thus, $T_{3r} = T_{t3r}$, $T_{2r} = T_{t2r}$, $T_{1r} = T_{t1r}$. For the static pressures $P_{2r}$ and $P_{3r}$, the law of perfect gas is used. 

\section{Optimization of the system}
\label{sec:optim}

\subsection{Formulation of the optimization problem}
\label{sec:formulation}

We consider an optimization problem using the 18 design variables given in Table~\ref{tab:variables}. All other design parameters and physical properties are fixed (see Table~\ref{tab:param} in appendix). Under the model developed in Section~\ref{sec:analysis}, the ECS is thus ruled by a system of 13 equations (Eq.~(1)--(13)) with 13 unknowns, which are the stagnation temperatures and pressures, the powers exchanged between the fluid and the compressor and turbine, and the rotational speed of the rotating ensemble: $T_{t2}, T_{t3}, T_{t4}, T_{t5}, T_{t2r}, T_{t3r}, P_{t2}, P_{t3}, P_{t4}, P_{t5}, \dot{W}_C, \dot{W}_T \text{\ and\ } \omega$. Eq.(14) and Eq.(15) give respectively the mass and the entropy generation rate, which are the objectives of the optimization. Additionally, we formulate the following 15 inequality constraints (see Table~\ref{tab:static}):
\begin{equation*}
\begin{array}{lcccccc}
c_{1-2} &:& T_{min} &\leq& T_5 &\leq& T_{max}, \\
c_{3-4} &:& P_{min} &\leq& P_5 &\leq& P_{max}, \\
c_{5-6} &:& 0.5 &\leq& \epsilon_1 &\leq& 0.9, \\
c_{7-8} &:& 0.5 &\leq& \epsilon_2 &\leq& 0.9, \\
c_{9} &:& \|C_2\| &\leq& 0.95 \sqrt{\gamma R T_2}, \\
c_{10} &:& \|C_3\| &\leq& 0.95 \sqrt{\gamma R T_3}, \\
c_{11} &:& \|C_4\| &\leq& 0.95 \sqrt{\gamma R T_4}, \\
c_{12} &:& \|C_5\| &\leq& 0.95 \sqrt{\gamma R T_5}, \\
c_{13} &:& \rc \omega &\leq& \sqrt{\gamma R T_3}, \\
c_{14} &:& \rt \omega &\leq& \sqrt{\gamma R T_4}, \\
c_{15} &:& \mathcal{P_{HT}} &\leq& \m c_p \left( T_c - T_5 \right).
\end{array}
\end{equation*}

The constraints $c_1$ to $c_4$ are standard specifications. The air injected into the cabin must lie between $T_{min} = 15^\circ \text{\ C}$ and $T_{max} = 25^\circ \text{\ C}$ and at a pressure close to the atmospheric pressure. Thus we take $P_{min} = 101.3 \text{\ kPa}$ and $P_{max} = 1.05P_{min}$. The constraints $c_5$ to $c_8$ are on the heat exchangers efficiencies. The design should be efficient enough but not too expensive to manufacture. The constraints $c_9$ to $c_{12}$ are on the air velocity. In the model, we have assumed that the air density remains constant throughout the bleed stream, which is inaccurate if the flow becomes supersonic. We take a 5\% margin to account for the possible variations of uncertain parameters and avoid numerical instabilities. Similarly, it is required via constraints $c_{13}$ and $c_{14}$ that the compressor and turbine blade tip speeds be subsonic. Constraint $c_{15}$ stems from the sizing scenario considered in this study: The dissipated power must be greater than the power produced by the passengers, the crew, the equipments and the environment (see Section~\ref{sec:scenario}). Note that an equality constraint is not necessary because the constraint is expected to be active at the optima.

To ensure the feasibility of the system and avoid numerical issues, we enforce the following restrictions on the design variables (see Figure~\ref{fig:parameters}):
\begin{equation*}
\begin{array}{lcccccc}
d_1 &:& \m &\leq& \mr , \\
d_2 &:& \bc &\leq& h_c \rc , \\
d_3 &:& \bt &\leq& h_t \rt , \\
d_{4-5} &:& \rpc + 0.02 &\leq& \rtc &\leq& \rc , \\
d_{6-7} &:& \rpt + 0.02 &\leq& \rtt &\leq& \rt , \\
d_8 &:& \Delta &\geq& 0 , \\
d_9 &:& \frac{\tan(\betac)}{\bc} &\geq& - \frac{\tan(\alphat)}{\bt}.
\end{array}
\end{equation*}
where $\Delta$ in $d_8$ is the discriminant of Eq~(7)--(9), seen as a second order polynomial equation in $\omega$. The conditions $d_8$ and $d_9$ are necessary to ensure that there exists a real solution $\omega > 0$ to Eq~(9). When two such solutions are possible, we take the largest one. Note that this parametrization implies that the optimization needs to be performed on a non-hypercubic design domain.

\begin{table*} 
\renewcommand{\arraystretch}{1.3}
\caption{Design variables description}
\label{tab:variables}
\begin{center}
\begin{tabular*}{\textwidth}{@{\extracolsep\fill}lcclcc@{}}
\hline
Description & Not. & Domain & Description & Not. & Domain\\
\hline
Bleed flowrate (kg.s\textsuperscript{-1}) & $\m$ & [2, 8] & Ram flowrate (kg.s\textsuperscript{-1}) & $\mr$ & [2, 8]\\
Compressor outlet radius (m) & $\rc$ & [0.1, 0.3] & Turbine inlet radius (m) & $\rt$ & [0.1, 0.3]\\
Compressor inlet foot radius (m) & $\rpc$ & [0.03, 0.1] & Turbine outlet foot radius (m) & $\rpt$ & [0.03, 0.1]\\
Compressor inlet tip radius (m) & $\rtc$ & [0.04, 0.2] & Turbine outlet tip radius (m) & $\rtt$ & [0.04, 0.2]\\
Compressor outlet blade height (m) & $\bc$ & [0.01, 0.1] & Turbine inlet blade height (m) & $\bt$ & [0.01, 0.1]\\
Compressor outlet angle (rad) & $\betac$ & [$-\frac{\pi}{3}$, $\frac{\pi}{3}$] & Turbine inlet angle (rad) & $\alphat$ & [0, $\frac{\pi}{3}$]\\
Heat exchanger 1: x length (m) & $\Lxa$ & [0.025, 0.7] & Heat exchanger 2: x length (m) & $\Lxb$ & [0.025, 0.7]\\
Heat exchanger 1: y length (m) & $\Lya$ & [0.025, 0.7] & Heat exchanger 2: y length (m) & $\Lyb$ & [0.025, 0.7]\\
Heat exchanger 1: z length (m) & $\Lza$ & [0.025, 0.7] & Heat exchanger 2: z length (m) & $\Lzb$ & [0.025, 0.7]\\
\hline
\end{tabular*}
\end{center}
\end{table*}

\subsection{Optimization algorithm}

The optimization is performed using the BMOO algorithm \cite{feliot2016bayesian}. This algorithm implements a Bayesian approach to the multi-objective optimization problem in the presence of non-linear constraints. The objectives and constraints of the problem are modeled using Gaussian process emulators (see, e.g., \cite{williams2006gaussian}) and the algorithm performs a sequential optimization procedure where the next sample is chosen as the maximizer of an extended version of the expected improvement sampling criterion (see, e.g., \cite{jones1998efficient}). In practice, this requires to solve an auxiliary optimization problem at each iteration. The BMOO algorithm uses sequential Monte Carlo techniques to conduct this auxiliary optimization (see, e.g., \cite{del2006sequential}). A population of candidate designs distributed according to a density of interest in the design space is sampled at each iteration and the maximizer of the extended expected improvement is chosen out of this population. 

We take advantage of this to handle non-hypercubic design domains by truncating the density of interest so as to propose only candidates that lie in the desired region. This is straightforward because sequential Monte Carlo methods do not require that the normalizing constant of the target density be known. For initialization, a pseudo maximin design of experiments on the non-hypercubic design domain (see e.g. \cite{auffray2012maximin}) can be achieved using rejection sampling. A large population of particles is sampled uniformly on the containing hypercube defined using the values of Table~\ref{tab:variables}. The particles which do not respect the constraints $d_1$ to $d_9$ are then discarded and the population of surviving particles is pruned until the desired population size is reached. During the pruning step, particles that are too close to other particles are discarded, thus raising the maximin distance. Note that in practice, this requires that the volume of the design domain be not too small compared with the volume of the containing hypercube (the ratio of volumes was estimated close to $6\%$ for this particular application).   

Because the computation of the objectives and constraints values for a given design requires to solve the non-linear system formed by Eq.(1) to Eq.(13), it may happen that no solution can be found, in which case it is not possible to provide values of the constraints and objectives for the design under study. Also, some designs can lead to supersonic solutions for which the values of temperatures and pressures predicted by the model can be inaccurate. When this happens, we prefer to consider such designs as simulation failures and not use the values returned by the model. In the optimization procedure, this is taken into account in order to prevent the optimizer to explore regions where simulation failures are likely, by multiplying both the sampling criterion of BMOO and the density in the sequential Monte Carlo procedure by a probability of observability. This technique has been proposed by Lee and co-authors \cite{lee2011optimization}. A statistical model is learned on the observed/non-observed data and provides a probability of satisfying the hidden constraints leading to simulation failures. In this work, a nearest-neighbours classifier using 5 neighbours and the $L_2$ distance is used to that purpose.

\subsection{Optimization results}

\begin{figure}
\centerline{
\includegraphics[height=7cm]{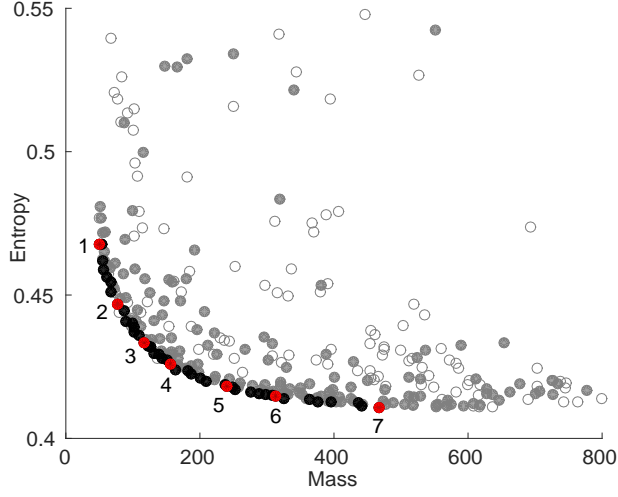}
}
\caption{Pareto front obtained with the BMOO optimizer using 500 samples. Empty circles are non-feasible solutions. Grey disks are feasible but dominated solutions. Black and red disks are feasible and non-dominated solutions.}
\label{fig:pareto}
\end{figure}
The algorithm is run with a limiting budget of $N_{max}=500$ calls to the simulation model, and an initial design of $N_{init} = 90$ samples. The set of optimal trade-off solutions found by the algorithm is shown on Figure~\ref{fig:pareto}. 

Among the initial design of experiments, 44 experiments led to simulation failures and 92 additional failures occurred during the optimization process. Further investigation revealed that most of the simulation failures occurred because the flow was supersonic in the compressor, which happens with high probability when the bleed flowrate is high and the compressor radii are low. Regarding the constraints satisfaction, no feasible observations were made in the initial sample and the algorithm found one after 25 iterations. 

The design parameters associated to 7 trade-off solutions chosen along the Pareto front are given in Table~\ref{tab:results}. Several observations can be made on these results. First we note that the bleed flowrate remains constant along the front. This is because $c_{15}$ is active (see Section~\ref{sec:formulation}) and $T_5 = T_{min}$ for optimal designs, which forces the value of the bleed flowrate. The ram flowrate is less constrained and varies along the front. We note its strong influence on the entropy generation rate (see Eq.(15)). The variation of the mass on the other hand mostly comes from the variations of $\Lza$ and $\Lzb$. As the heat exchangers height is raised, the entropy generation rate is lowered but the mass augments. The values of $\Lxa$, $\Lya$, $\Lxb$, and $\Lyb$ are set close to their maximal values, which permits to achieve efficiencies between 0.7 and 0.8. Note that the pressure losses are assumed constant in this study. Further work is required to better understand their impact on the entropy generation rate when the heat exchangers dimensions become small. Regarding the turbine and compressor dimensions, they are set as small as possible, which keeps the mass low and augments the fluid velocity, thus achieving good performances.

\begin{table*}[h!] 
\renewcommand{\arraystretch}{1.05}
\caption{Optimal design variables values found by the optimization algorithm for the points 1 to 7 (in red) of Figure~\ref{fig:pareto}. The values of the most influential variables are in bold.}
\label{tab:results}
\begin{center}
\begin{tabular*}{\textwidth}{@{\extracolsep\fill}lccccccc@{}}
\hline
 & 1 & 2 & 3 & 4 & 5 & 6 & 7 \\
\hline
$\m$ & 2.95 & 2.92 & 2.94 & 2.94 & 2.94 & 2.95 & 2.94 \\
$\mr$ & $\bm{7.74}$ & $\bm{6.86}$ & $\bm{5.63}$ & $\bm{5.06}$ & $\bm{4.64}$ & $\bm{4.40}$ & $\bm{4.27}$ \\
$\rpc$ & 0.07 & 0.05 & 0.05 & 0.03 & 0.03 & 0.07 & 0.04 \\
$\rtc$ & 0.10 & 0.08 & 0.08 & 0.08 & 0.06 & 0.09 & 0.10 \\
$\rc$ & 0.10 & 0.11 & 0.10 & 0.10 & 0.12 & 0.12 & 0.13 \\
$\bc$ & 0.01 & 0.01 & 0.05 & 0.05 & 0.04 & 0.02 & 0.03 \\
$\betac$ & 0.36 & 0.74 & 0.97 & -0.16 & 0.61 & 0.94 & 0.48 \\
$\rpt$ & 0.03 & 0.03 & 0.03 & 0.03 & 0.03 & 0.03 & 0.03 \\
$\rtt$ & 0.05 & 0.05 & 0.05 & 0.05 & 0.05 & 0.05 & 0.05 \\
$\rt$ & 0.10 & 0.10 & 0.11 & 0.12 & 0.11 & 0.10 & 0.11 \\
$\bt$ & 0.02 & 0.02 & 0.04 & 0.02 & 0.04 & 0.03 & 0.03 \\
$\alphat$ & 1.04 & 0.50 & 0.89 & 1.01 & 0.44 & 0.79 & 0.30 \\
$\Lxa$ & 0.67 & 0.65 & 0.68 & 0.68 & 0.63 & 0.69 & 0.70 \\
$\Lya$ & 0.65 & 0.68 & 0.61 & 0.67 & 0.67 & 0.66 & 0.65 \\
$\Lza$ & $\bm{0.03}$ & $\bm{0.04}$ & $\bm{0.07}$ & $\bm{0.12}$ & $\bm{0.17}$ & $\bm{0.20}$ & $\bm{0.32}$ \\
$\Lxb$ & 0.66 & 0.69 & 0.66 & 0.66 & 0.70 & 0.68 & 0.69 \\
$\Lyb$ & 0.69 & 0.53 & 0.68 & 0.65 & 0.65 & 0.68 & 0.65 \\
$\Lzb$ & $\bm{0.03}$ & $\bm{0.06}$ & $\bm{0.09}$ & $\bm{0.10}$ & $\bm{0.17}$ & $\bm{0.25}$ & $\bm{0.36}$ \\
\hline
$\mathcal{M}$ & 49.78 & 77.13 & 117.00 & 156.57 & 240.03 & 312.40 & 466.69 \\
$\dot{\mathcal{S}}$ & 0.47 & 0.45 & 0.43 & 0.43 & 0.42 & 0.41 & 0.41 \\
\hline

\end{tabular*}
\end{center}
\end{table*}

\section{Conclusions}
\label{sec:conclusion}

In this article, a one dimensional model of the environment control system of a commercial aircraft is proposed. The model permits to emulate the behaviour of the system when the geometries of its components vary, for a scenario where the aircraft is on ground, full of passengers, equipments running, and with an outside temperature of $50^\circ$C. The system is optimized using the BMOO algorithm, which implements a Bayesian approach to the multi-objective optimization problem in the presence of non-linear constraints, and trade-off design solutions in terms of mass and entropy generation rate of the system are identified.

As a particularity, the optimization is performed on a non-hypercubic design domain and involves hidden constraints. This is a situation that is often encountered in engineering design optimization. The BMOO algorithm is successfully adapted to this new setup, which makes it possible to conduct a multi-objective optimization using a reasonable number of calls to the numerical simulation model.

The BMOO algorithm is primarily designed to address problems where the computational time associated to the model evaluation is high, which is not the case here. In this study, most of the computational time required to conduct the optimization was taken by the optimizer and more work is needed to make the implementation of the algorithm more efficient. Nevertheless, the algorithm achieves very satisfactory results and is a competitive algorithm to address multi-objective optimization problems with several constraints.

~\\

\textbf{Acknowledgements:} This research work has been carried out within the Technological Research Institute SystemX, using public funds from the French Programme \emph{Investissements d'Avenir}. We also thank Airbus Group Innovations for their contributions to this work.

\bibliographystyle{plain}
\bibliography{biblio}

\begin{thebibliography}{10}

\bibitem{auffray2012maximin}
Y.~Auffray, P.~Barbillon, and J.~M. Marin.
\newblock Maximin design on non hypercube domains and kernel interpolation.
\newblock {\em Statistics and Computing}, 22(3):703--712, 2012.

\bibitem{bejan2001need}
A.~Bejan and D.~L. Siems.
\newblock The need for exergy analysis and thermodynamic optimization in
  aircraft development.
\newblock {\em Exergy, An International Journal}, 1(1):14--24, 2001.

\bibitem{del2006sequential}
P.~Del~Moral, A.~Doucet, and A.~Jasra.
\newblock Sequential {M}onte {C}arlo samplers.
\newblock {\em Journal of the Royal Statistical Society: Series B (Statistical
  Methodology)}, 68(3):411--436, 2006.

\bibitem{feliot2016bayesian}
P.~Feliot, J.~Bect, and E.~Vazquez.
\newblock A {B}ayesian approach to constrained single- and multi-objective
  optimization.
\newblock {\em Journal of Global Optimization}, 37 pages. 2016. In press.
  Available online at http://dx.doi.org/10.1007/s10898-016-0427-3.

\bibitem{fonseca1998multiobjective}
C.~M. Fonseca and P.~J. Fleming.
\newblock Multiobjective optimization and multiple constraint handling with
  evolutionary algorithms. {I}. {A} unified formulation.
\newblock {\em IEEE Transactions on Systems, Man and Cybernetics. Part A:
  Systems and Humans}, 28(1):26--37, 1998.

\bibitem{jones1998efficient}
D.~R. Jones, M.~Schonlau, and W.~J. Welch.
\newblock Efficient global optimization of expensive black-box functions.
\newblock {\em Journal of Global Optimization}, 13(4):455--492, 1998.

\bibitem{lee2011optimization}
H.~K.~H. Lee, R.~B. Gramacy, C.~Linkletter, and G.~A. Gray.
\newblock Optimization subject to hidden constraints via statistical emulation.
\newblock {\em Pacific Journal of Optimization}, 7(3):467--478, 2011.

\bibitem{perez2002optimization}
I.~P{\'e}rez-Grande and T.~J. Leo.
\newblock Optimization of a commercial aircraft environmental control system.
\newblock {\em Applied thermal engineering}, 22(17):1885--1904, 2002.

\bibitem{williams2006gaussian}
C.~Rasmussen and C.~K.~I. Williams.
\newblock {G}aussian processes for machine learning.
\newblock {\em MIT Press}, 2(3):1--266, 2006.

\bibitem{vargas2001thermodynamic}
J.~V.~C. Vargas and A.~Bejan.
\newblock Thermodynamic optimization of finned crossflow heat exchangers for
  aircraft environmental control systems.
\newblock {\em International Journal of Heat and Fluid Flow}, 22(6):657--665,
  2001.

\end{thebibliography}

\pagebreak
\appendix

\section{Additional material}

\begin{table}[h!] 
\renewcommand{\arraystretch}{1.5}
\caption{Static properties equations}
\label{tab:static}
\begin{center}
\begin{tabular*}{\textwidth}{@{\extracolsep\fill}lll@{}}
\hline
Static temperatures & Static pressures & Fluid velocities\\
\hline
$T_5 = T_{t5} - \frac{C_{5x}^2}{2}$ & 
$\frac{P_{t5}}{P_5} = \left( 1 + \frac{\gamma-1}{2} \left(\frac{C_{5x}}{\gamma R T_5}\right)^2 \right)^{\frac{\gamma}{\gamma-1}}$ & 
$C_{5x} = \frac{\m}{\pi \left( \rtt^2-\rpt^2 \right) \rho}$ \\

$T_4 = T_{t4} - \frac{C_{4m}^2}{2}$ & 
$\frac{P_{t4}}{P_4} = \left( 1 + \frac{\gamma-1}{2} \left(\frac{C_{4m}}{\gamma R T_4}\right)^2 \right)^{\frac{\gamma}{\gamma-1}}$ & 
$C_{4m} = \frac{\m}{2\pi \rt \bt \rho \cos \alphat}$ \\

$T_3 = T_{t3} - \frac{C_{3m}^2}{2}$ & 
$\frac{P_{t3}}{P_3} = \left( 1 + \frac{\gamma-1}{2} \left(\frac{C_{3m}}{\gamma R T_3}\right)^2 \right)^{\frac{\gamma}{\gamma-1}}$ & 
$C_{3m} = \sqrt{\left( \frac{\rc \omega - \m \tan \betac}{2\pi \rc \bc \rho}\right)^2 + \left( \frac{\m}{2\pi \rc \bc \rho} \right)^2}$ \\

$T_2 = T_{t2} - \frac{C_{2x}^2}{2}$ & 
$\frac{P_{t2}}{P_2} = \left( 1 + \frac{\gamma-1}{2} \left(\frac{C_{2x}}{\gamma R T_2}\right)^2 \right)^{\frac{\gamma}{\gamma-1}}$ & 
$C_{2x} = \frac{\m}{\pi \left( \rtc^2-\rpc^2 \right) \rho}$ \\
\hline
\end{tabular*}
\end{center}
\end{table}
\begin{table*}[h!] 
\renewcommand{\arraystretch}{1.25}
\caption{Parameters values used in the experiments of Section~\ref{sec:optim}}
\label{tab:param}
\begin{center}
\begin{tabular*}{\textwidth}{@{\extracolsep\fill}lcclcc@{}}
\hline
Description & Not. & Value & Description & Not. & Value \\
\hline

\multicolumn{6}{l}{\hspace{0.5cm}\emph{Simulation parameters}} \\
Ambient temperature (K) & $T_a$ & 323 & Ambient pressure (Pa) & $P_a$ & 101.3e3 \\
Number of passengers & $N_{pax}$ & 120 & Number of crewmembers & $N_{crew}$ & 5 \\
Thermal power passengers (W) & $\mathcal{P}_{pax}$ & 70 & Thermal power crew (W) & $\mathcal{P}_{crew}$ & 100 \\
Thermal power equipments (W) & $\mathcal{P}_{eq}$ & 4800 & Outside flow dissipation (W)& $\mathcal{P}_{out}$ & 3000 \\
Bleed temperature (K) & $T_1$ & 473 & Bleed pressure (Pa) & $P_1$ & 260e3 \\
Pressure losses (Pa) & $\Delta P_{HX}$ & 40e3 & Valve opening & $\theta$ & 0 \\
Ram stream cross surface (m\textsuperscript{2}) & $A_r$ & 0.20 & Fan efficiency & $\eta_F$ & 0.95 \\
Air specific heat (J.kg\textsuperscript{-1}.K\textsuperscript{-1}) & $c_p$ & 1004 & Air isentropic coefficient & $\gamma$ & 1.4 \\
Perfect gaz constant (J.kg\textsuperscript{-1}.K\textsuperscript{-1}) & $R$ & 287 & & & \\

\multicolumn{6}{l}{\hspace{0.5cm}\emph{Heat exchangers}} \\
Viscosity bleed (kg.m\textsuperscript{-1}.s\textsuperscript{-1}) & $\mu$ & 2.28e-5 & Viscosity ram (kg.m\textsuperscript{-1}.s\textsuperscript{-1}) & $\mu_r$ & 2.28e-5 \\
H.T. ratio bleed stream (m\textsuperscript{-1}) & $\beta$ & 2231 & H.T. ratio ram stream (m\textsuperscript{-1}) & $\beta_r$ & 1115 \\
Plate spacing bleed stream (m) & $b$ & 5.21e-3 & Plate spacing ram stream (m) & $b_r$ & 12.3e-3 \\
Prandtl number bleed stream & $Pr$ & 0.7 & Prandtl number ram stream & $Pr_r$ & 0.7 \\
Hydraulic diameter bleed (m) & $Dh$ & 1.54e-3 & Hydraulic diameter ram (m) & $Dh_r$ & 3.41e-3 \\ 
Convection length bleed (m) & $\lambda$ & 0.035 & Convection length ram (m) & $\lambda_r$ & 0.035 \\
Representative density (kg.m\textsuperscript{-3}) & $\rho_{HX}$ & 1415 & Fin thickness (m) & $\delta$ & 0.102e-3 \\
Wall thickness (m) & $t_w$ & 6e-4 & Thermal conductivity (W.m\textsuperscript{-1}.K\textsuperscript{-1}) & $k_w$ & 237 \\

\multicolumn{3}{l}{\hspace{0.5cm}\emph{Compressor}} & \multicolumn{3}{l}{\hspace{0.5cm}\emph{Turbine}} \\
Adiabatic efficiency & $\eta_c$ & 0.8 & Adiabatic efficiency & $\eta_t$ & 0.92 \\
Aspect ratio & $h_c$ & 0.7 & Aspect ratio & $h_t$ & 0.5 \\
Blades thickness (m) & $e_c$ & 0.01 & Blades thickness (m) & $e_t$ & 0.01 \\
Number of blades & $Z_c$ & 21 & Number of blades & $Z_t$ & 21 \\
\hline
\end{tabular*}
\end{center}
\end{table*}

\end{document}